\documentclass[12pt,reqno]{article}

\usepackage[usenames]{color}
\usepackage{amssymb}
\usepackage{graphicx}
\usepackage{amscd}

\usepackage[colorlinks=true,
linkcolor=webgreen,
filecolor=webbrown,
citecolor=webgreen]{hyperref}

\definecolor{webgreen}{rgb}{0,.5,0}
\definecolor{webbrown}{rgb}{.6,0,0}

\usepackage{color}
\usepackage{fullpage}
\usepackage{float}

\usepackage{graphics,amsmath,amssymb}
\usepackage{amsthm}
\usepackage{amsfonts}
\usepackage{latexsym}
\usepackage{epsf}

\newcommand{\seqnum}[1]{\href{http://www.research.att.com/cgi-bin/access.cgi/as/~njas/sequences/eisA.cgi?Anum=#1}{\underline{#1}}}

\makeatletter
\def\rdots{\mathinner{\mkern1mu\raise\p@
\vbox{\kern7\p@\hbox{.}}\mkern2mu
\raise4\p@\hbox{.}\mkern2mu\raise7\p@\hbox{.}\mkern1mu}}
\makeatother

\begin{document}

\begin{center}
\epsfxsize=4in
\end{center}

\begin{center}
\vskip 1cm{\LARGE\bf
COMPOSITION\\
\vskip .1in
OF ORDINARY GENERATING\\
\vskip .1in
FUNCTIONS}
\vskip 1cm
\large
Vladimir Kruchinin\\
Tomsk State University\\
of Control Systems and Radioelectronics\\
Tomsk\\
Russion Federation\\
\href{mailto:kru@ie.tusur.ru}{\tt kru@ie.tusur.ru}\\
\end{center}

\vskip .2in

\begin{abstract}
A solution is proposed for the problem of composition of ordinary generating functions. A new class of functions that provides a composition of ordinary generating functions is introduced; main theorems are presented; compositae are written for polynomials, trigonometric and hyperbolic functions, exponential and log functions. It is shown that the composition holds true for many integer sequences.
\end{abstract}

\newtheorem{theorem}{Theorem}
\newtheorem{corollary}[theorem]{Corollary}
\newtheorem{lemma}[theorem]{Lemma}
\newtheorem{proposition}[theorem]{Proposition}
\newtheorem{conjecture}[theorem]{Conjecture}
\newtheorem{defin}[theorem]{Definition}
\newenvironment{definition}{\begin{defin}\normalfont\quad}{\end{defin}}
\newtheorem{examp}[theorem]{Example}
\newenvironment{example}{\begin{examp}\normalfont\quad}{\end{examp}}
\newtheorem{rema}[theorem]{Remark}
\newenvironment{remark}{\begin{rema}\normalfont\quad}{\end{rema}}

\abstractname{ A solution is proposed for the problem of composition of ordinary generating functions. A new class of functions that provides a composition of ordinary generating functions is introduced; main theorems are presented; compositae are written for polynomials, trigonometric and hyperbolic functions, exponential and log functions. It is shown that the composition holds true for many integer sequences}

\section{Introduction}

Generating functions are an efficient tool of solving mathematical problems. Given the ordinary generating functions $F(x)=\sum\limits_{n\geqslant 1} f(n)x^n$ and $R(x)=\sum\limits_{n\geqslant 0} r(n)x^n$, the operation of composition of generating functions $A(x)=R(F(x))$ is defined correctly. \cite{stanley-vol2,Marc,ConcreteMath,Egor}. However, coefficients of the composition of generating functions are difficult to find. Stanley \cite{stanley-vol2} came close to the solution of the problem and proposed a formula for the composition of exponential generating functions based on ordered partitions of a finite set. 
	Let us show that the basis for the composition of ordinary generating functions is ordered partitions of a positive integer $n$ and put forward basic formulae for the coefficients of the composition of ordinary generating functions. For this purpose, we introduce several definitions.
	
\begin{defin} An ordinary generating function $F(x)$ is a series that belongs to the ring of formal power series in one variable $K[[x]]$:
$$
F(x)=\sum_{n\geq 0}f(n)x^n,
$$
where $f(n):P\rightarrow K$, $P$ is a set of nonnegative numbers; $K$ is a commutative field. 
\end{defin}

Further we consider only ordinary generating functions. The known generating functions are denoted as $F(x)$, $R(x)$, $G(x)$, and the desired generating function as $A(x)$.

\begin{defin} An ordered partition (composition) of a positive integer $n$ is an ordered sequence of positive integers $\lambda_i$ such that
$$
\sum_{i=1}^k\lambda_i=n, 
$$
where $\lambda_i\geq1$ and $k=\overline{1,n}$ are parts of the ordered partition. 
\end{defin} 
\noindent ${C_n}$ is a set of all ordered partitions of $n$.\\ 
$\pi_k\in C_n$ is an ordered partition of $C_n$ with $k$ parts.\\
The ordered partitions of $n$ have been much studied \cite{Andrews,Kru2010}. 

\section{Compositae and their properties}

Let there be functions $f(n)$ and $r(n)$ and their generating functions $F(x)=\sum_{n\geqslant 1} f(n)x^n$, $R(x)=\sum_{n\geqslant 0} r(n)x^n$. Then, calculating the composition of the generating functions $A(x)=R(F(x))$ requires \cite{Egor}
\begin{equation}\label{KruForm1}
[F(x)]^k=\sum_{n\geq k}\sum_{\lambda_i>0 \atop {\lambda_1+\lambda_2+\ldots+\lambda_k=n}}
f(\lambda_1)f(\lambda_2)\ldots f(\lambda_k)x^n.
\end{equation}
Hence it follows that for the function $a(n)$ of the composition of generating functions with $n>0$, the formula 
$$
a(0)=r(0),
$$
\begin{equation}\label{KruForm2}
a(n)=\sum_{k=1}^n\left[\sum_{\lambda_i>0 \atop {\lambda_1+\lambda_2+\ldots+\lambda_k=n}}
f(\lambda_1)f(\lambda_2)\ldots f(\lambda_k)\right]r(k)
\end{equation}	
holds true. Further the composition of generating functions is written implying that $a(0)=r(0)$. 

\begin{remark} It should be  noted that the summation in formulae (\ref{KruForm1}),(\ref{KruForm2}) is over all ordered partitions of $n$ that have exactly $k$ parts, because $\{\lambda_1+\lambda_2+\ldots+\lambda_k=n\},~\lambda_i>0,~ i=\overline{1,k}$ (further we use the reduction $\pi_k\in C_n$). 
\end{remark}

Thus, the ordered partitions of $n$ are the basis for calculation of the composition of generating functions.

Let us consider the following example. Assume that $f(0)=0$, $f(n)=1$ for all $n>0$. This function is defined by the generating function $F(x)=\frac{x}{1-x}$. Then, the expression
$$
\sum_{\pi_k\in C_n}
f(\lambda_1)f(\lambda_2)\ldots f(\lambda_k)	   
$$
gives the number of ordered partitions of $n$ with exactly $k$ parts; this number is equal to ${n-1\choose k-1}$ \cite{Andrews}. Thus,
$$
\sum_{\pi_k \in C_n}f(\lambda_1)f(\lambda_2)\ldots f(\lambda_k)={n-1\choose k-1}.	   
$$
Hence it follows that the formula valid for any generating function $R(x)=\sum_{n\geqslant 0} r(n)x^n$ and $A(x)=R\left(\frac{x}{1-x}\right)$ is  
$$
a(n)=\sum_{k=1}^n{n-1\choose k-1}r(k).
$$
\begin{examp}
For $R(x)=\frac{x}{1-x}$, we have the composition $A(x)=\frac{x}{1-2x}$ and
$$
a(n)=\sum_{k=1}^n{n-1\choose k-1}=2^{n-1}.
$$
Thus, we calculate the total number of ordered partitions of $n$.
\end{examp}

\begin{examp} We have $R(x)=e^x$, then for the composition $A(x)=e^{\frac{x}{1-x}}$ we can write
$$
a(n)=\sum_{k=1}^n{n-1\choose k-1}\frac{1}{k!}
$$
(see A000262 formula Herbert S. Wilf).
\end{examp}

\begin{examp} We have $R(x)=\frac{x}{1-x-x^2}$, then for the composition $A(x)=R(\frac{x}{1-x})$ we can write
$$
a(n)=\sum_{k=1}^n{n-1\choose k-1}F(k),
$$
where $F(k)$ is the Fibonacci numbers (see A001519, formula Benoit Cloitre). 
\end{examp}

\begin{defin} A composita\index{Composita! Definition} of the generating function $F(x)=\sum_{n>0}f(n)x^n$ is the function 
\begin{equation}
\label{Fnk0}F^{\Delta}(n,k)=\sum_{\pi_k \in C_n}{f(\lambda_1)f(\lambda_2)\ldots f(\lambda_k)}.
\end{equation}
\end{defin}

Calculation of $F^{\Delta}(n,k)$ is of prime importance to obtain a composition of generating functions, because from formula (\ref{KruForm2}) it follows that the formula valid for the composition $A(x)=R(F(x))$ is
\begin{equation}
\label{Superposita}
a(n)=\sum_{k=1}^n F^{\Delta}(n,k) r(k).
\end{equation} 
The basis for the derivation of a composita is calculation of the ordered partition $\pi_k$ of $C_n$. From formula (\ref{KruForm1}) it follows that the generating function of the composita is equal to
$$
[F(x)]^k=\sum_{n\geq k}F^{\Delta}(n,k)x^n.
$$
For $F(x)$, the condition $f(0)=0$ holds true, and hence numbering for the composita begins with $k=1$, $n=1$. For $k=1$, $F^{\Delta}(n,k)=f(n)$. For $k>n$, $F^{\Delta}(n,k)$ is equal to zero. This statement stems from the fact that there is no ordered partition of $n$ in which the number of parts is larger than $n$. 

The above example demonstrates that the Pascal triangle is a composita for the generating function $\frac{x}{1-x}$\index{Composita of a generating function! $F(x)=\frac{x}{1-x}$} and deriving the composition $A(x)=R\left(\frac{x}{1-x}\right)$ requires the use of
$$
F^{\Delta}(n,k)={n-1 \choose k-1}.
$$

Let us derive a recurrence formula for the composita of a generating function.

\begin{theorem} \label{TheoremFnk}
For the composita $F^{\Delta}(n,k)$ of the generating function $F(x)=\sum_{n>0}f(n)x^n$, the following relation holds true:
\begin{equation}
\label{Composita}
F^{\Delta}(n,k)=\left\{
\begin{array}{ll}
f(n),& k=1,\\
{[f(1)]}^n,& k=n,\\
\sum_{i=0}^{n-k}f(i+1)F^{\Delta}(n-i-1,k-1)& k<n. \\
\end{array}
\right.
\end{equation}
\end{theorem}

\begin{proof}  

Let us derive a recurrence formula for the $c_{n,k}$ number of ordered partitions of $n$ that have exactly $k$ parts. Let us introduce the operation $pos[\lambda^*,\pi_k]$ of adjunction of the new part $\lambda^*$ on the left to a certain ordered partition $\pi_k \in C_n$ providing that $\lambda^*>0$. From the ordered partition $\pi_k \in C_n$ this operation obtains an ordered partition $\pi_{k+1} \in C_{\lambda^*+n}$. Let us extend this operation to sets. Assume that $C_{n,k}=\{\pi_k|\pi_k\in C_n\}$, then the set $\hat C = pos[\lambda^*,C_{n,k}]$ is a subset $C_{\lambda^*+n,k+1}$. Thus, we can write   
$$
C_{n,k}=pos[1,C_{n-1,k-1}]\cup pos[2,C_{n-1,k-1}]\cup \ldots\cup pos[n-k-1,C_{k-1,k-1}].
$$   
In this case, the condition
$$
pos[i,C_{n-i,k-1}]\cap pos[j,C_{n-j,k-1}]=\oslash
$$
is fulfilled for all $i\neq j$, because the first parts of the ordered partitions $\pi_k$ do not coincide. Hence,
\begin{equation}
\label{cnk}
c_{n,k}=\sum_{i=0}^{n-k}c_{n-i-1,k-1},
\end{equation}
and $c_{k,k}=1$ because we have the only ordered partition $\pi_k=\{1+1+\ldots+1=n\}$, and $c_{n,1}=1$ because $\pi_1=\{n=n\}$.

Let us now consider expression (\ref{Fnk0}). Using expression (\ref{cnk}), we can write
$$
\begin{array}{l}
F^{\Delta}(n,k)=\\
=f(1)F^{\Delta}(n-1,k-1)+f(2)F^{\Delta}(n-2,k-1)+\ldots+\\
+f(n-k+1)F^{\Delta}(k-1,k-1).
\end{array}
$$
The set $C_{n,n}$ consists of the only ordered partition $\{1+1+\ldots +1\}$, and then $F_{n,n}^{\Delta}=[f(1)]^n$; the set $C_{n,1}$ consists of $\{n\}$, and then $F_{n,1}^{\Delta}=f(n)$.
Thus, the theorem is proved.
\end{proof}

Consideration of formula (\ref{Superposita}) allows the conclusion that the composita does not depend on $R(x)$ and characterizes the generating function $F(x)$. In tabular form, the composita is represented as
$$
\begin{array}{ccccccccccc}
&&&&& F_{1,1}^{\Delta}\\
&&&& F_{2,1}^{\Delta} && F_{2,2}^{\Delta}\\
&&& F_{3,1}^{\Delta} && F_{3,2}^{\Delta} && F_{3,3}^{\Delta}\\
&& F_{4,1}^{\Delta} && F_{4,2}^{\Delta} && F_{4,3}^{\Delta} && F_{4,4}^{\Delta}\\
& \rdots && \vdots && \vdots && \vdots && \ddots\\
F_{n,1}^{\Delta} && F_{n,2}^{\Delta} && \ldots && \ldots && F_{n,n-1}^{\Delta} && F_{n,n}^{\Delta}\\
\end{array}
$$
or, knowing that $F_{1,n}^{\Delta}=f(n)$, $F_{n,n}^{\Delta}={[f(1)]}^n$, as 
$$
\begin{array}{ccccccccccc}
&&&&& f(1)\\
&&&& f(2) && f^2(1)\\
&&& f(3) && F_{3,2}^{\Delta} && f^3(1)\\
&& f(4) && F_{4,2}^{\Delta} && F_{4,3}^{\Delta} && f^4(1)\\
& \rdots && \vdots && \vdots && \vdots && \ddots\\
f(n) && F_{n,2}^{\Delta} && \ldots && \ldots && F_{n,n-1}^{\Delta} && f^n(1)\\
\end{array}
$$ 

Below are the terms of the composite of the generating function $F(x)=\frac{x}{1-x}$:
$$
\begin{array}{ccccccccccc}
&&&&&1\\
&&&& 1 && 1\\
&&& 1 && 2 && 1\\
&& 1 && 3 && 3 && 1\\
& 1 && 4 && 6 && 4 && 1\\
1 && 5 && 10 && 10 && 5 && 1
\end{array}
$$

\begin{theorem}
For a given ordinary generating function $F(x)=\sum_{n\geq 1}{f(n)x^n}$, the composita $F^{\Delta}(n,k)$ always exists and is unique.
\end{theorem}
\begin{proof} Without proof. \end{proof}

\section{Calculation of compositae}

Calculation of compositae is based on derivation of the generating function of a composita
$$
[A(x)]^k=\sum_{n\geqslant k} A^{\Delta}(n,k)x^n
$$
and operation on them.

\begin{theorem} \label{theorem_ak} Let there be a generating function $F(x)=\sum_{n>0} f(n)x^n$, its composita $F^{\Delta}(n,k)$, and a constant $\alpha$. Then, the generating function $A(x)=\alpha F(x)$ has the composita  
$$
A^{\Delta}(n,k)=\alpha^k F^{\Delta}(n,k).
$$
\end{theorem}
\begin{proof}
$$
[A(x)]^k=[\alpha F(x)]^k=\alpha^k[F(x)]^k.
$$
\end{proof}

\begin{theorem} \label{theorem_an} Let there be a generating function $F(x)=\sum_{n>0} f(n)x^n$, its composita $F^{\Delta}(n,k)$, and a constant $\alpha$. Then, the generating function $A(x)=F(\alpha x)$ has the composita  
$$
A^{\Delta}(n,k)=\alpha^n F^{\Delta}(n,k).
$$
\end{theorem}
\begin{proof}
By definition, we have
$$
A^{\Delta}(n,k)=\sum_{\pi_k \in C_n}{\alpha^{\lambda_1}f(\lambda_1)\alpha^{\lambda_2}f(\lambda_2)\ldots \alpha^{\lambda_k}f(\lambda_k)}=
$$
$$
=\alpha^{n}\sum_{\pi_k \in C_n}{f(\lambda_1)f(\lambda_2)\ldots f(\lambda_k)}=\alpha^{n}F^{\Delta}(n,k).
$$
\end{proof}

\begin{theorem} \label{theorem_mult} Let there be a generating function $F(x)=\sum_{n>0} f(n)x^n$, its composita $F^{\Delta}(n,k)$, a generating function $B(x)=\sum_{n\geqslant 0} b(n)x^n$ and $[B(x)^k]=\sum_{n\geqslant 0}B(n,k)x^n$. Then, the generating function $A(x)=F(x)B(x)$ has the composita
$$
A^{\Delta}(n,k)=\sum_{i=k}^{n} F^{\Delta}(i,k)B(n-i,k).
$$
\end{theorem}
\begin{proof}
 Because $a(0)=f(0)b(0)=0$, $A(x)$ has the composita $A^{\Delta}(n,k)$. On the other hand,
$$
 [A(x)]^k=[F(x)]^k[B(x)]^k.
$$ 
This, reasoning from the rule of product of generating functions, gives
$$
A^{\Delta}(n,k)=\sum_{i=k}^{n} F^{\Delta}(i,k)B(n-i,k).
$$
\end{proof}
For $B(x)$  $b(0)=0$, the formula has the form:
$$
A^{\Delta}(n,k)=\sum_{i=k}^{n-k} F^{\Delta}(i,k)B^{\Delta}(n-i,k).
$$

\begin{theorem} \label{theorem_sum} Let there be generating functions $F(x)=\sum_{n>0} f(n)x^n$, $G(x)=\sum_{n>0} g(n)x^n$  and their compositae $F^{\Delta}(n,k)$ , $G^{\Delta}(n,k)$. Then, the generating function $A(x)=F(x)+G(x)$ has the composita
$$
A^{\Delta}(n,k)=F^{\Delta}(n,k)+\sum_{j=1}^{k-1}{k\choose j}\sum_{i=j}^{n-k+j}F^{\Delta}(i,j)G^{\Delta}(n-i,k-j)+G^{\Delta}(n,k).
$$
\end{theorem}
\begin{proof}
According to the binomial theorem, we have
$$
[A(x)]^k=\sum_{j=0}^k{{k\choose j} [F(x)]^j[G(x)]^{k-j}},
$$
$$
[F(x)]^j=\sum_{n\geqslant j}F^{\Delta}(n,j),
$$
and 
$$
[G(x)]^{k-j}=\sum_{n\geqslant k-j}G^{\Delta}(n,k-j).
$$
According to the rule of multiplication of series, we obtain
$$
A^{\Delta}(n,k)=F^{\Delta}(n,k)+\sum_{j=1}^{k-1}{k\choose j}\sum_{i=j}^{n-k+j}F^{\Delta}(i,j)G^{\Delta}(n-i,k-j)+G^{\Delta}(n,k).
$$
\end{proof}

\begin{defin} Let there be a composition of generating functions $A(x)=R(F(x))$. Then, the product of two compositae will be termed a composite of the composition $A(x)$ and denoted as
$A^{\Delta}(n,k)=F^{\Delta}(n,k)\circ R^{\Delta}(n,k)$.
\end{defin}

\begin{theorem}\label{CompozitProduct}  Let there be two generating functions $F(x)=\sum_{n>0}f(n)x^n$ and $R(x)=\sum_{n>0}r(n)x^n$, and their compositae $F^{\Delta}(n,k)$ and $R^{\Delta}(n,k)$. Then, the expression valid for the product of the compositae $A^{\Delta}=F^{\Delta}\circ R^{\Delta}$ is
\begin{equation} \label{CompositaARF}
A^{\Delta}(n,m)=\sum_{k=m}^n F^{\Delta}(n,k)R^{\Delta}(k,m).
\end{equation} 
\end{theorem}

\begin{proof}
$$
[A(x)]^m=[G(F(x)]^m=G^m(F(x)).
$$
Hence, according to the composition rule and taking into account that the nonzero terms $G^{\Delta}(n,m)$ begin with $n\geqslant m$, we have
$$
A^{\Delta}(n,m)=\sum_{k=m}^{n}F^{\Delta}(n,k)G^{\Delta}(k,m).
$$
\end{proof}
{\bf Corollary }. Because the composition of generating functions is an associative operation and
$$
F(x)\circ(R(x)\circ G(x))=(F(x)\circ R(x))\circ G(x),
$$
the product of compositae is also an associative operation and 
\begin{equation}
\sum_{k=m}^n \sum\limits_{i=k}^n F^{\Delta}(n,i)R^{\Delta}(i,k)G^{\Delta}(k,m)=\nonumber
\sum_{k=m}^n \sum\limits_{i=k}^n R^{\Delta}(n,i)G^{\Delta}(i,k)F^{\Delta}(k,m).\nonumber 
\end{equation}
\section{Compositae of generating functions}

\subsection{Identical composita} \index{Composita! Identical}

\begin{defin} An identical composita $Id^{\Delta}(n,k)$ is a composita of the generating function $F(x)=x$. 
\end{defin}
By definition, $[F(x)]^k=x^k$. Then
\begin{equation}
F^{\Delta}(n,k)=\left\{ 
\begin{array}{ll}
1,& n=k,\\
0,& n\neq k. 
\end{array}
\right.
\end{equation}
Thus, $F^{\Delta}(n,k)=\delta_{n,k}$, where $\delta_{n,k}$ is the Kronecker delta.
It is easily seen that for any generating function $A(x)$, we have the identity
$$
a(n)=\sum_{k=1}^nId^{\Delta}(n,k)a(k).
$$
The composita of the function $F(x)=x^m$ is
\begin{equation}
F^{\Delta}(n,k,m)=\delta_{\frac{n}{m},k}, \mod(n,m)=0~\hbox{or}~n=km.
\end{equation}

\subsection{Compositae of polynomials}

\subsubsection{Composita for $P_2(x)=(ax+bx^2)$ }  

Let us consider $P_2(x)=(ax+bx^2)$. Then, $p_2(0)=0$,  $p_2(1)=a$ and $p_2(2)=b$, and the rest are $p_2(n)=0,~n>2$. The composita of the function $F(x)=a x$ is equal to $a^k\delta_{n,k}$, and the composita of the function $G(x)=bx^2$ is equal to $b^k\delta_{\frac{n}{2},k}$. Using sum theorem (\ref{theorem_sum}), we obtain
$$
P_2^{\Delta}(n,k)=\sum_{j=0}^k {k \choose j}\sum_{i=j}^{n-k+j} a^j\delta_{i,j}b^{k-j}\delta_{\frac{n-i}{2},k-j},
$$
$\delta_{\frac{n-i}{2},k-j}=1$ for  $\frac{n-i}{2}=k-j$, whence $i=n-2k+2j$. So we have
$$
P_2^{\Delta}(n,k)=\sum_{j=0}^k {k \choose j} a^j\delta_{n-2k+2j,j}b^{k-j}.
$$
Now $\delta_{n-2k+2j,j}=1$ for $n-2k+2j=j$, whence $j=2k-n$. So we obtain
\begin{equation} \label{com_alphabeta}
P_2^{\Delta}(n,k,a,b)={k \choose n-k} a^{2k-n}b^{n-k}
\end{equation}
for $\lceil\frac{n}{2}\rceil\leq k \leq n$.

Thus, the composition $A(x)=R(ax+bx^2)$ can be found using the expression:
$$
a(n)=\sum_{k=\lceil\frac{n}{2}\rceil}^n{k \choose n-k}a^{2k-n}b^{n-k}r(k).
$$
For example, let us derive an expression for the coefficients of the generating function $A(x)=e^{x+\frac{1}{2}x^2}$
(see \seqnum{A000085} ). Taking into account that this function is an exponential generating function, we obtain 
$$
a(n)=n!\sum_{k=\lceil\frac{n}{2}\rceil}^n{k \choose n-k}\frac{1}{2^{n-k}}\frac{1}{k!}.
$$
Another example is $A(x)=R(F(x))$, where $R(x)=\frac{x}{1-x}$ and $F(x)=x+x^2$, $A(x)=\frac{x+x^2}{1-x-x^2}$. Hence
$$
a(n)=\sum_{k=\lceil\frac{n}{2}\rceil}^n{k \choose n-k}
$$
(see \seqnum{A000045}).
 
\subsubsection{Composita for $P_3(x)=ax+bx^2+cx^3$  }  

The polynomial $P_3(x)=ax+bx^2+cx^3$ can be expressed as
$$
P_3(x)=ax+xP_2(x,b,c).
$$
The composita $ax$ is equal to $\delta(n,k)a^k$, and the composita $xP_2(x)$ to $A_2{\Delta}(n-k,k)$. Then, on the strength of the theorem on the composita of the sum of generating functions, we have 
$$
A_3^{\Delta}(n,k,a,b,c)=\sum_{j=0}^k{k\choose j}\sum_{i=j}^{n-k+j}A_2{\Delta}(i-j,j,b,c)\delta(n-i,k-j)a^{k-j}.
$$
Simplification gives $\delta(n-i,k-j)=1$ for $n-i=k-j$, whence we have $i=n-k+j$ and
$$
A_3^{\Delta}(n,k,a,b,c)=\sum_{j=0}^k{k\choose j}A_2(n-k,j,b,c)a^{k-j},
$$ 
where $A_2^{\Delta}(n-k,j,b,c)={j \choose n-k-j} b^{2j+k-n}b^{n-k-j}$. Hence,
$$
A_3^{\Delta}(n,k,a,b,c)=\sum_{j=0}^k{k\choose j}{j \choose n-k-j}a^{k-j}b^{2j+k-n}b^{n-k-j}.
$$ 
Then, for the generating function $A(x)=\frac{1}{1-ax-bx^2-cx^3}$, the following expression holds true:
$$
a(n)=\sum_{k=1}^n\sum_{j=0}^k{k\choose j}{j \choose n-k-j}a^{k-j}b^{2j+k-n}b^{n-k-j}.
$$
 
\subsubsection{Composita for $P(x)=ax+cx^3$} 

An important condition in the foregoing examples is that $a,b,c\neq 0$. Therefore, if $b=0$ the formula for the composita should be sought for over again. For example,
$$
P(x)=ax+cx^3.
$$
In this case, the expression for the composita is
$$
P^{\Delta}(n,k)={k\choose \frac{3k-n}{2}}a^{\frac{3k-n}{2}}c^{\frac{n-k}{2}},
$$
where $(n-k)$ is exactly divisible by 2.
For example, for the generating function $A(x)=\frac{1}{1-x-x^3}$ the following expression holds true:
$$
a(n)=\sum_{k=1}^n {k\choose \frac{3k-n}{2}}
$$
(see \seqnum{A000930}).

\subsubsection{Composita for $P_4(x)=ax+bx^2+cx^3+dx^4$} 

At $n=4$, the polynomial $P_4(x)=ax+bx^2+cx^3+dx^4$ can be expressed as
$$
P_4(x)=P_2(x)+x^2P_2(x).
$$
The generating function of the composita for $x^2P_2(x)$ is equal to
$$
x^{2k}{k \choose n-k}c^{2k-n}b^{n-k}x^n={k \choose n-k}c^{2k-n}b^{n-k}x^{n+2k},
$$
and hence the expression for the composita is
$$
{k \choose n-3k}c^{4k-n}b^{n-3k}.
$$
Then the composita $P_4(x)$ has the following expression:
{
$$
\sum_{j=0}^{k}{k\choose j}\sum_{i=j}^{n-k+j}{j\choose i-j}a^{2j-i}b^{i-j}{k-j \choose n-i-3(k-j)}c^{4(k-j)-(n-i)}d^{n-i-k+j}.
$$
For example, for the generating function $A(x)=\frac{1}{1-ax-bx^2-cx^3-dx^4}$ the following expression holds true:
$$
a(n)=\sum_{k=1}^n \sum_{j=0}^{k}{k\choose j}\sum_{i=j}^{n-k+j}{j\choose i-j}a^{2j-i}b^{i-j}{k-j \choose n-i-3(k-j)}c^{4(k-j)-(n-i)}d^{n-i-k+j}.
$$
At $a=b=c=d=1$, we obtain the generating function $A(x)=\frac{1}{1-x-x^2-x^3-x^4}$. Hence
$$
a(n)=\sum_{k=1}^n \sum_{j=0}^{k}{k\choose j}\sum_{i=j}^{n-k+j}{j\choose i-j}{k-j \choose n-i-3(k-j)}.
$$
\subsubsection{Composita for $P_5(x)=ax+bx^2+cx^3+dx^4+ex^5$} 

For finding the composita of the $m$th power polynomial, we can propose the recurrent algorithm
$$
A_m^{\Delta}(n,k)=\sum_{j=0}^k A_{m-1}^{\Delta}(n-k,j)a^{k-j},  
$$
providing that $A_{m-1}(n,0)=1$. Using this recurrent algorithm, we obtain the composita of the 5th power polynomial:
$$\sum_{r=0}^{k}{a^{k-r}\,{{k}\choose{r}}\,\sum_{m=0}^{r}{b^{r-m}\,
 \left(\sum_{j=0}^{m}{c^{m-j}\,d^{r-n+m+k+2\,j}\,e^{v}\,{{j
 }\choose{v}}\,{{m}\choose{j}}}\right)\,{{r}\choose{m}}}}$$,
where $v=-r+n-m-k-j$.

\subsection{Composita for $A(x)=(\frac{ax}{1-bx})$ } 
 \label{com_pascalalphabeta}
 
For the generating function $F(x)=\frac{x}{(1-x)}$, $F^{\Delta}(n,k)={n-1 \choose k-1}$, and
$$ 
A(x)=(ab^{-1}\frac{bx}{1-bx}).
$$ Using theorems (\ref{theorem_ak},\ref{theorem_an}), we obtain
$$
A^{\Delta}(n,k)={n-1 \choose k-1} a^{k}b^{n-k}.
$$

\subsection{Compositae of the exponent}

Let us find the expression for the coefficients of the generating function $[B(x)]^k=e^{kx}$:
$$
B(x)^k=e^{xk}=\sum_{n\geqslant 0} \frac{k^n}{n!},
$$
whence it follows that
$$
B(n,k)=\frac{k^n}{n!}.
$$
Now, for $A(x)=xe^x$ the composita is equal to
\begin{equation}\label{Composita_xexp}
A^{\Delta}(n,k)=B(n-k,k)=\frac{k^{n-k}}{(n-k)!}.
\end{equation}

Let us write the composita for the generating function $A(x)=e^x-1$:
$$
A(x)^k=\sum_{m=0}^k {k \choose m} e^{mx}(-1)^{k-m},
$$
whence it follows that the composita is
\begin{equation}\label{Composita_expm1}
A^{\Delta}(n,k)=\sum_{m=0}^k {k \choose m} \frac{m^n}{n!}(-1)^{k-m}=\frac{k!}{n!}S_2(n,k),
\end{equation}
where $S_2(n,k)$ is the Stirling numbers of the second kind.
For the generating functions of the Bell numbers $A(x)=e^{e^x-1}$, we have 
$$
a(n)=n!\sum_{k=1}^n S_2(n,k)\frac{k!}{n!}\frac{1}{k!}=\sum_{k=1}^nS_2(n,k)
$$
(see \seqnum{A000110}).

\subsection{Composita for $\ln(1+x)$}
\label{Compozita_ln} 
\index{Composita of the generating function! $F(x)=\ln(x+1)$ }

Let $F(x)=\ln(x+1)$. Then, knowing the relation \cite{ConcreteMath} 
$$
\sum_{n=k}^\infty S_1(n,k) \frac{x^n}{n!} = \frac{\left[\ln(1+x)\right]^k}{k!},
$$
where $S_1(n,k)$ is the Stirling numbers of the first kind, and using formula (\ref{KruForm1}), we obtain the expression for the composita of the generating function $\ln(1+x)$:
\begin{equation}
F^{\Delta}(n,k)=\frac{k!}{n!}S_1(n,k).
\end{equation}

\subsection{Composita for the generating function of the Bernoulli numbers}

The generating function of the Bernoulli numbers is
$$
A(x)=\frac{x}{e^x-1}.
$$
This function can be represented as the composition $B(F(x))$, where $B(x)=\frac{\ln x}{x}$, $F(x)=e^x-1$. Let us find the expression for the coefficients of the generating function $[B(x)]^k$:
$$
[B(x)]^k=\sum_{n\geqslant 0} S_1(n,k)\frac{k!}{n!}x^{n-k},
$$
whence 
$$
B(n,k)=S_1(n+k,k)\frac{k!}{(n+k)!}.
$$
Knowing the composita of the function $F(x)$ (see \ref{Composita_expm1}),

$$
F^{\Delta}(n,k)=\frac{k!}{n!}S_2(n,k).
$$

Let us write the composition of the generating functions $A(x)^k=[B(e^x-1)]^k$: 
$$
A(n,m)=\left\{
\begin{array}{ll}
1, & n=0,\\
\sum_{k=1}^n S_2(n,k)\frac{k!}{n!} S_1(k+m,m)\frac{m!}{(k+m)!},& n>0.
\end{array}
\right.
$$
Then the composita of $xA(x)$ is
$$
A^{\Delta}(n,m)=\left\{
\begin{array}{ll}
1, & n=m,\\
\frac{m!}{(n-m)!}\sum_{k=1}^{n-m} \frac{k!}{(k+m)!}S_1(k+m,m)S_2(n-m,k),& n>m.
\end{array}
\right.
$$

\subsection{Composita for the generating function of the Fibonacci numbers}

Let us find the composita for the generating function of the Fibonacci numbers:
$$
A(x)=\frac{x}{1-x-x^2}.
$$
The function can be represented as the composition of the generating functions $A(x)=R(F(x))$, where $R(x)=\frac{x}{1-x}$, $F(x)=\frac{x}{1-x^2}$.
Let us find the composita for $F(x)$:
$$
F^{\Delta}(n,k)=\left\{
\begin{array}{ll}
{\frac{n+k}{2}-1\choose k-1}, & \hbox{at $n+k$ -- even},\\
0,& \hbox{at $n+k$ -- odd}.\\
\end{array}
\right.
$$
Now, using the operation of product of compositae, we find the composita of the generating function $A(x)$:
$$
A^{\Delta}(n,m)=\begin{array}{ll}
\sum\limits_{k=m}^n{\frac{n+k}{2}-1\choose k-1}{k-1\choose m-1}, & \hbox{at $n+k$ -- even}.\\
\end{array}
$$
Below are the first terms of the composita for the generating function of the Fibonacci numbers:
$$
\begin{array}{cccccccccccccc}
&&&&&&1\\
&&&&&1&&1\\
&&&&2&&2&&1\\
&&&3&&5&&3&&1\\
&&5&&10&&9&&4&&1\\
&8&&20&&22&&14&&5&&1\\
13&&38&&51&&40&&20&&6&&1\\
\end{array}
$$

\subsection{Composita for the generalized Fibonacci numbers}

Let us find the composita of the generating function: 
$$
F(x)=x+x^2+\ldots+x^m=\frac{x-x^{m+1}}{1-x}.
$$
Let us write $F(x)$ as the product of the functions $G(x)=x-x^{m+1}$ and $R(x)=\frac{1}{1-x}$.
Let us find the composita for $G(x)$. For this purpose, we consider the compositae of the functions $y(x)=x$ and $z(x)=-x^m$.
For $y(x)$, the composita is equal to $Id(n,k)=\delta_{n,k}$. For $z(x)=-x^m$, the composita is
$$
Z^{\Delta}(n,k)=(-1)^k\delta_{\frac{n}{m},k}.
$$
Then, on the strength of the theorem on the composite of sum of generating functions $y(x)+z(x)$, we have
$$
G^{\Delta}(n,k)=\sum_{j=0}{k \choose j}\sum_{i=j}^{n-k+j}Id(i,j)Z^{\Delta}(n-i,k-j)=
$$
$$
=\sum_{j=0}{k \choose j}\sum_{i=j}^{n-k+j}\delta_{i,j}\delta_{\frac{n-i}{m},k-j}(-1)^{k-j}.
$$
The function $\delta_{i,j}=1$ is only for $i=j$, and hence
$$
G^{\Delta}(n,k)=\sum_{j=0}{k \choose j}\delta_{\frac{n-j}{m},k-j}(-1)^{k-j}.
$$
The function $\delta_{\frac{n-j}{m},k-j}=1$ is only for $\frac{n-j}{m}=k-j$, and hence
$$
G^{\Delta}(n,k)={k\choose \frac{(m+1)k-n}{m}}(-1)^{\frac{n-k}{m}}.
$$
It is known that $R(n,k)={n+k-1\choose k-1}$. Then, with regard to the rule of finding the composita of the product of generating functions (case 2), we obtain
$$
F^{\Delta}(n,k)=\sum_{i=k}^n {k\choose \frac{(m+1)k-i}{m}}(-1)^{\frac{i-k}{m}}{n-i+k-1\choose k-1}.
$$
Let us consider the composita of the generating functions: 
$$
A(x)=\frac{F(x)}{1-F(x)}=\frac{x-x^{m+1}}{1-2x-x^{m+1}}.
$$ 
Hence, using the theorem on the product of compositae, we obtain the composita of the generating function $A(x)$:
\begin{eqnarray}
A^{\Delta}(n,l)=\sum_{k=l}^n F^{\Delta}(n,k){k-1 \choose m-1}=\nonumber\\
=\sum_{k=m}^n \sum_{i=k}^n {k\choose \frac{(m+1)k-i}{m}}(-1)^{\frac{i-k}{m}}{n-i+k-1\choose k-1}{k-1 \choose l-1}\nonumber.
\end{eqnarray}
For $l=1$, we derive the formula for the generalized Fibonacci numbers:
\begin{equation}
F_n^{(m)}=\sum_{k=1}^n \sum_{i=k}^n {k\choose \frac{(m+1)k-i}{m}}(-1)^{\frac{i-k}{m}}{n-i+k-1\choose k-1}.
\end{equation}

\subsection{Composita of the generating function for the Catalan numbers} 

Let $F(x)=x\frac{1-\sqrt{1-4x}}{2x}$, then the composita has the form
$$
F^{\Delta}(n,k)=\sum_{i=0}^{n-k} C(i)F_{n-i-1,k-1}^{\Delta},
$$ where $C(i)$ is the Catalan numbers. The composita $F^{\Delta}(n,k)$ has the following triangular form:
$$
\begin{array}{ccccccccccc}
&&&&&1\\
&&&& 1 && 1\\
&&& 2 && 1 && 1\\
&& 5 && 5 && 3 && 1\\
& 14 && 14 && 9 && 4 && 1\\
\end{array}
$$
Let us consider the sequence \seqnum{A009766} called the Catalan triangle. This triangle is given by the formula
$$
a(n,m)={n+m\choose n}\frac{n-k+1}{n+1}.
$$
Below are the initial values of the triangle, and $n$ and $m$ begin with zero.  
$$
\begin{array}{ccccccccccc}
&&&&&1\\
&&&& 1 && 1\\
&&& 2 && 1 && 2\\
&& 1 && 3 && 5 && 5\\
& 1 && 4 && 9 && 14 && 14\\
\end{array}
$$
Comparison of two triangles suggests that $a(n,k)=F^{\Delta}(n+1,n-k+1)$. Hence, the composita for the Catalan generating function is equal to
$$
F^{\Delta}(n,k)={2n-k-1\choose n-1}\frac{k}{n}.
$$

Thus, the expression valid for the coefficients of the composition $A(x)=R(\frac{1-\sqrt{1-4x}}{2})$ is
$$
a(n)={2n-k-1\choose n-1}\frac{k}{n}\cdot r(k).
$$

\subsection{Composita of the generating function $\frac{x}{\sqrt{1-x}}$}

This generating function can be represented as the composition of the functions: 
$$
\frac{x}{\sqrt{1-x}}=x\frac{1}{1-\left(2\frac{\sqrt{1-\frac{4x}{4}}}{2}-1\right)}=x\frac{1}{1-2C(\frac{1}{4}x)},
$$ 
where $C(x)=\frac{1-\sqrt{1-4x}}{2}$.

Using the formula of composition, we finally obtain
$$
A^{\Delta}(n,m)=\left\{
\begin{array}{ll}
1, & n=m,\\
\sum_{k=1}^{n-m}{2n-2m-k-1\choose n-m-1}\frac{k}{n-m}2^{k-2n+2m}{k+m-1 \choose m-1},& n>m.
\end{array}
\right.
$$

\subsection{Compositae of trigonometric functions}

\subsubsection{Composita of the sine}

Using the expression
$$
\sin(x)=\frac{e^{ix}-e^{-ix}}{2i},
$$
we obtain $\sin(x)^k$:
$$
\sin(x)^k=\frac{1}{2^ki^k}\sum_{m=0}^k {k \choose m} e^{imx}e^{-i(k-m)x}(-1)^{k-m}=\frac{1}{2^ki^k}\sum_{m=0}^k {k \choose m} e^{i(2m-k)x}(-1)^{k-m}.
$$
Hence the composita is
$$
\frac{1}{2^k}i^{n-k}\sum_{m=0}^k {k \choose m} \frac{(2m-k)^n}{n!}(-1)^{k-m}.
$$

Taking into account that $n-k$ is an even number and the function is symmetric about $k$, we obtain the composita of the generating function $\sin (x)$:
$$
A^{\Delta}(n,k)=\left\{
\begin{array}{ll}
\frac{1}{2^{k-1}n!}\sum_{m=0}^{\lfloor \frac{k}{2} \rfloor} {k \choose m} (2m-k)^n(-1)^{\frac{n-k}{2}-m}, & (n-k) -\hbox{even}\\
0, & (n-k) - \hbox{odd}
\end{array}
\right.
$$

\begin{examp} For the Euler numbers we know the exponential generating function $\frac{1}{1-sin(x)}$.
Hence,
$$
E_{n+1}=\sum_{{k=1} \atop \hbox{\scriptsize $n+k$ even} }^n \frac{1}{2^{k-1}}\sum_{m=0}^{\lfloor \frac{k}{2} \rfloor} {k \choose m} (2m-k)^n(-1)^{\frac{n+k}{2}-m}
$$
(see \seqnum{A000111}).
\end{examp}

\begin{examp} For the generating function $A(x)=e^{\sin(x)}$, the valid expression is
$$
a_n=\sum_{{k=1} \atop \hbox{\scriptsize $n+k$ even} }^n \frac{1}{2^{k-1}k!}\sum_{m=0}^{\lfloor \frac{k}{2} \rfloor} {k \choose m} (2m-k)^n(-1)^{\frac{n+k}{2}-m}
$$
(see \seqnum{A002017}).
\end{examp}

\subsubsection{Compositae of the cosine}

Knowing that
$$
\cos(x)=\frac{e^{ix}+e^{-ix}}{2},
$$
We have
$$
[\cos x]^k=\frac{1}{2^k}\sum_{j=0}^k {k \choose j} e^{(2j-k)ix}=
$$
$$
=\frac{1}{2^k}\sum_{n\geqslant 0}\sum_{j=0}^k {k \choose j}(2j-k)^ni^n\frac{x^n}{n!}. 
$$
Hence
$$
B(n,k)=\frac{1}{2^kn!}(-1)^{\frac{n}{2}}\sum_{j=0}^k {k \choose j}(2j-k)^n.
$$
Then, the composita of the generating function $xcos(x)$ is
$$
A^{\Delta}(n,k)=\left\{
\begin{array}{ll}
\frac{1}{2^{k-1}(n-k)!}(-1)^{\frac{n-k}{2}}\sum_{j=0}^k {k \choose j}(2j-k)^{n-k}, & n-k~-\hbox{even}\\
0, & n-k~-~\hbox{odd}.
\end{array}
\right.
$$
The composita of the function $cos(x)-1$ is equal to
$$
A^{\Delta}(n,k)=\sum_{i=0}^k B(n,i)(-1)^{k-i}.
$$

Let us consider the following example. Let there be a generating function $A(x)=sec(x)=\frac{1}{cos(x)}=\frac{1}{1+(\cos(x)-1)}$. Hence, on the strength of the formula of composition and composita $(\cos(x)-1)$, we obtain
$$a(n)=\sum_{k=1}^{2\,n}{\sum_{m=0}^{k}{{{k}\choose{m}}\,2^{1-m}\,\left(
 \sum_{j=0}^{{{m}\over{2}}}{\left(2\,j-m\right)^{2\,n}\,{{m}\choose{j
 }}}\right)\,\left(-1\right)^{n+m}}}
$$
(see \seqnum{A000364}).

\subsubsection{Composita for $\tan(x)$}

For the tangent, we know the identity
$$
\tan(x)=\frac{e^{ix}-e^{-ix}}{i(e^{ix}-e^{-ix})}.
$$
Division of the numerator and denominator by $e^{ix}$ gives
$$
\tan(x)=\frac{1-e^{-2ix}}{i(1-e^{-2ix})}.
$$
Multiplication of the numerator and denominator by i, and addition and then subtraction of unity gives
$$
\tan(x)=i\frac{e^{-2ix}-1}{2-(e^{-2ix}-1)}.
$$
Whence it follows that
$$
\tan(x)=\frac{i}{2}\frac{e^{-2ix}-1}{1-\frac{1}{2}(e^{-2ix}-1)}.
$$
Thus, the function $\tan(x)$ is expressed as the composition of the functions
$$
F(x)=\frac{i}{2}\frac{x}{1+\frac{1}{2}x}
$$
and functions $R(x)=e^{-2ix}-1$.
The composita for $F(x)$ is equal to
$$
F^{\Delta}(n,k)=\frac{1}{2^n}(-1)^{n-k}{n-1 \choose k-1}i^k.
$$
The composita for $R(x)$ is equal to
$$
R^{\Delta}(n,k)=(-2i)^n\frac{k!}{n!}S_2(n,k),
$$
where $S_2(n,k)$ is the Stirling numbers of the second kind.
Then, on the strength of the theorem on the product of compositae, we obtain the composita of the function $\tan(x)$:
$$
A^{\Delta}(n,m)=\sum_{k=m}^n (-2i)^nS_2(n,k)\frac{k!}{n!}\frac{1}{2^k}(-1)^{k-m}{k-1 \choose m-1}i^m.
$$
After transformation, we obtain
$$
A^{\Delta}(n,m)=(-1)^{\frac{n+m}{2}}\sum_{k=m}^n (2)^{n-k}S_2(n,k)\frac{k!}{n!}(-1)^{n+k-m}{k-1 \choose m-1}.
$$
Then at $k=1$, the expression for the tangential numbers is
$$
a(n)=\left(-1\right)^{n+1}\,\sum_{j=1}^{2\,n+1}{\left(-1\right)^{j}\,j!
 \,2^{2\,n-j+1}\,S_2(2\,n+1,j)}
$$
(see \seqnum{A000182})

Let us consider the example $A(x)=e^{\tan(x)}$:
$$
a(n)=\sum_{k=1}^{n}{{{\left(-1\right)^{{{n+k}\over{2}}}\,\sum_{j=k}^{n}{
 {{j-1}\choose{k-1}}\,j!\,2^{n-j}\,\left(-1\right)^{n-k+j}\,
 {\it S_2}\left(n , j\right)}}\over{k!}}}
$$
(see \seqnum{A006229}).
For more examples, see \seqnum{A000828},\seqnum{A000831},\seqnum{A003707}

\subsubsection{Composita for $x^2\cot(x)$}

It is known that 
$$
x^2\cot(x)=ix\frac{e^{ix}+e^{-ix}}{e^{ix}-e^{-ix}}=ix^2+\frac{2ix^2}{e^{2ix}-1}.
$$
The composita $ix^2$ is equal to $\delta(\frac{n}{2},k)i^k$, and the composita for $\frac{2ix^2}{e^{2ix}-1}$ is equal to
$$
(2i)^{n-k}B^{\Delta}(n,k),
$$
where $B^{\Delta}(n,k)$ is the composita for the generating function of the Bernoulli numbers.
Using the theorem on the composita of the sum of generating functions, we obtain the composita of the function $x^2\cot(x)$:
$$
A^{\Delta}(n,k)=\delta(\frac{n}{2},k)i^k+\sum_{j=1}^k  B^{\Delta}(n-2k+2j,j)(2i)^{n-2k+j}i^{k-j}=
$$
$$
=\delta(\frac{n}{2},k)i^k+i^{n-k}\sum_{j=1}^k  B^{\Delta}(n-2k+2j,j)2^{n-2k+j}
$$

\subsubsection{Composita of the arc tangent $F(x)=\arctan(x)$}
Let us consider the generating function of the arc tangent:
$$
\arctan(x)=\sum_{n\geq 0}\frac{(-1)^n}{(2n+1)}x^{2n+1}. 
$$
Let us find an expression for the composita of the arc tangent from the operation of product of compositae. For this purpose, the expression
$$
\arctan(x)=\frac{i}{2}(\ln(1-ix))-\ln(1+ix))
$$
is written as follows:
$$
\arctan(x)=\frac{i}{2}\ln(1-\frac{2ix}{1+ix}).
$$
The composita of the function $f(x)=\frac{2ix}{1+ix}$ is equal to
$$
F^{\Delta}(n,k)=2^k{n-1 \choose k-1}i^n,
$$
whence it follows that
\begin{equation}
A_z^{\Delta}(n,m)=\frac{i^m}{2^m}\sum_{k=m}^n 2^k{n-1 \choose k-1}i^n\frac{m!}{k!}S_1(k,m).
\end{equation}

\begin{equation}
A_z^{\Delta}(n,m)=\frac{(-1)^{\frac{m+n}{2}}}{2^m}\sum_{k=m}^n 2^k{n-1 \choose k-1}\frac{m!}{k!}S_1(k,m).
\end{equation}

Below are the first terms of the composita of the arc tangent $A^{\Delta}(n,k)$ in the triangular form:
$$
\begin{array}{cccccccccccccccc}
&&&&&&&1\\
&&&&&&0&&1\\
&&&&&-\frac{1}{3}&&0&&1\\
&&&&0&&-\frac{2}{3}&&0&&1\\
&&&\frac{1}{5}&&0&&-1&&0&&1\\
&&0&&\frac{23}{45}&&0&&-\frac{4}{3}&&0&&1\\
&-\frac{1}{7}&&0&&\frac{14}{15}&&0&&-\frac{5}{3}&&0&&1\\
\end{array}
$$

\begin{examp} Let there be $R(x)=\frac{1}{1-x}$, then the coefficients of the generating function   
$$
A(x)=\frac{1}{1-\arctan(x)}
$$
are expressed by the formula:
$$
a(n)=\sum_{m=1}^n \frac{(-1)^{\frac{m+n}{2}}}{2^m}\sum_{k=m}^n 2^k{n-1 \choose k-1}\frac{m!}{k!}S_1(k,m).
$$
Hence, summation of rows of the composita of the arc tangent gives the following series:
$$
A(x)=1+x+x^2+\frac{2}{3}x^3+\frac{1}{3}x^4+\frac{1}{5}x^5+\frac{8}{45}x^6+\ldots.
$$
\end{examp}

\begin{examp} Let $A(x)=e^{\arctan(x)}$, then the valid expression is
$$
a(n)=n!\,\sum_{m=1}^{n}{{{\left(-1\right)^{{{3\,n+m}\over{2}}}\,\sum_{i=
 m}^{n}{{{2^{i}\,S_1(i,m)\,{{n-1}\choose{i-
 1}}}\over{i!}}}}\over{2^{m}}}}
$$
(see \seqnum{A002019}).
\end{examp}

\subsection{Compositae of hyperbolic functions}

For the hyperbolic sine, we have the known expression:
$$
\sinh(x)=\frac{e^x-e^{-x}}{2}.
$$

Let us find the composita of this generating function. For this purpose, we write
$$
\left(\frac{e^x-e^{-x}}{2}\right)^k=\frac{1}{2^k}\left(e^x+e^{-x}\right)^k=\frac{1}{2^k}\sum_{i=0}^k {k\choose i} e^{(k-i)x}(-1)^ie^{-ix}=
$$
$$
=\frac{1}{2^k}\sum_{i=0}^k (-1)^i{k\choose i} e^{(k-2i)x}.
$$
Let us write $e^x$ as a series, then we obtain
$$
\frac{1}{2^k}\sum_{i=0}^k (-1)^i{k\choose i} \sum_{n\geqslant 0}\frac{(k-2i)^n}{n!}x^n.
$$
Hence, the composita is \index{Composita of the generating function! $F(x)=\sinh(x)$}
$$
F^{\Delta}(n,k)=\frac{1}{2^k}\sum_{i=0}^k (-1)^i{k\choose i} \frac{(k-2i)^n}{n!}.
$$
For example, for $A(x)=e^{\sinh x}$ the valid expression is
$$
a(n)=\sum_{k=1}^{n}{{{\sum_{i=0}^{k}{\left(-1\right)^{i}\,\left(k-2\,i
 \right)^{n}\,{{k}\choose{i}}}}\over{2^{k}\,k!}}}
$$
(see \seqnum{A002724}).

For the hyperbolic cosine, we have
$$
\cosh(x)=\frac{e^x+e^{-x}}{2}.
$$
Then, 
\begin{eqnarray}
\cosh^k(x)=\left(\frac{e^x+e^{-x}}{2}\right)^k=\frac{1}{2^k}\sum_{i=0}^k {k\choose i} e^{(k-2i)x}=\nonumber\\
=\frac{1}{2^k}\sum_{i=0}^k {k\choose i} \sum_{n\geqslant 0}\frac{(k-2i)^n}{n!}x^n\nonumber,
\end{eqnarray}
and hence the composita of the generating function $x\cosh(x)$ is
$$
F^{\Delta}(n,k)=\frac{1}{2^k}\sum_{i=0}^k {k\choose i} \frac{(k-2i)^{n-k}}{(n-k)!}.
$$
For example, for $A(x)=e^{\cosh x}$ the valid expression is
$$
\sum_{k=1}^{n}{{{\left(\sum_{i=0}^{k}{\left(k-2\,i\right)^{n-k}\,{{
k}\choose{i}}}\right)\,{{n}\choose{k}}}\over{2^{k}}}}
$$
see \seqnum{A003727}).

\section{Conclusion}

The operation of the composition $A(x)=R(F(x))$ of ordinary generating functions requires:

1. Finding the composita $F^{\Delta}(n,k)$ of the generating function $F(x)$ with the use of theorems (\ref{theorem_ak},\ref{theorem_an},\ref{theorem_mult},\ref{theorem_sum},\ref{CompozitProduct})

2. Writing the composition in the form
$$
a(n)=\sum_{k=1}^n F^{\Delta}(n,k)r(n).
$$

\bigskip
\hrule
\bigskip

\noindent 2000 {\it Mathematics Subject Classification}:
Primary 05A15; Secondary 30B10.

\noindent \emph{Keywords: } Composition of ordinary generating function, ordered partitions, composita of ordinary generating function, integer sequence.

\bigskip
\hrule
\bigskip

\noindent (Concerned with sequences
\seqnum{A000045},\
\seqnum{A000085},
\seqnum{A000110},
\seqnum{A000111},
\seqnum{A000182},
\seqnum{A000262},
\seqnum{A000364},
\seqnum{A000828},
\seqnum{A000831},
\seqnum{A000930},
\seqnum{A001519},
\seqnum{A002017},
\seqnum{A002019},
\seqnum{A002714},
\seqnum{A003707},
\seqnum{A003727},
\seqnum{A006229},
\seqnum{A009766}.)

\bigskip
\hrule
\bigskip

\end{document}